\documentclass[11pt]{article}
\usepackage{fullpage}

\def\stocmode{0}
\def\jamesmode{0}
\def\arxivmode{0}
\def\fastmode{0}

\def\showauthornotes{0}

\def\showkeys{0}
\def\showdraftbox{1}
\def\showcolorlinks{1}
\def\usemicrotype{1}
\def\showfixme{1}

\ifnum\fastmode=0
\usepackage{etex}
\fi

\ifnum\fastmode=0
\usepackage[l2tabu, orthodox]{nag}
\fi

\usepackage{xspace,enumerate}

\usepackage[dvipsnames]{xcolor}

\ifnum\fastmode=0
\usepackage[T1]{fontenc}
\usepackage[full]{textcomp}
\fi

\usepackage[american]{babel}

\usepackage{mathtools}

\usepackage{amsthm}

\newtheorem{theorem}{Theorem}[section]
\newtheorem*{theorem*}{Theorem}

\newtheorem*{proposition*}{Proposition}
\newtheorem{lemma}[theorem]{Lemma}
\newtheorem*{lemma*}{Lemma}
\newtheorem{corollary}[theorem]{Corollary}
\newtheorem*{conjecture*}{Conjecture}

\newtheorem*{fact*}{Fact}

\newtheorem*{exercise*}{Exercise}

\newtheorem*{hypothesis*}{Hypothesis}

\theoremstyle{definition}

\newtheorem{exercise-easy}[theorem]{Exercise}
\newtheorem{exercise-med}[theorem]{Exercise}
\newtheorem{exercise-hard}[theorem]{Exercise$^\star$}

\newtheorem*{claim*}{Claim}

\newtheorem{remark}[theorem]{Remark}
\newtheorem*{remark*}{Remark}

\newtheorem*{observation*}{Observation}

\ifnum\arxivmode=0
\usepackage{newpxtext} %
\usepackage{textcomp} %
\usepackage[varg,bigdelims]{newpxmath}
\usepackage[scr=rsfso]{mathalfa}%
\usepackage{bm} %
\let\mathbb\varmathbb
\fi

\ifnum\arxivmode=1
\usepackage[varg]{pxfonts} %
\usepackage{textcomp} %
\usepackage[scr=rsfso]{mathalfa}%
\usepackage{bm} %
\fi

\ifnum\showkeys=1
\usepackage[color]{showkeys}
\fi

\makeatletter
\@ifpackageloaded{hyperref}
{}
{\usepackage[pagebackref]{hyperref}}

\definecolor{bleudefrance}{rgb}{0.01, 0.1, 1.0}
\definecolor{azure}{rgb}{0.0, 0.5, 1.0}

\ifnum\showcolorlinks=1
\hypersetup{
pagebackref=true,
colorlinks=true,
urlcolor=blue,
linkcolor=blue,
citecolor=OliveGreen}
\fi

\ifnum\showcolorlinks=0
\hypersetup{
pagebackref=true,
colorlinks=false,
pdfborder={0 0 0}}
\fi

\usepackage{prettyref}

\newcommand{\savehyperref}[2]{\texorpdfstring{\hyperref[#1]{#2}}{#2}}

\newrefformat{eq}{\savehyperref{#1}{\textup{(\ref*{#1})}}}
\newrefformat{lem}{\savehyperref{#1}{Lemma~\ref*{#1}}}
\newrefformat{def}{\savehyperref{#1}{Definition~\ref*{#1}}}
\newrefformat{thm}{\savehyperref{#1}{Theorem~\ref*{#1}}}
\newrefformat{cor}{\savehyperref{#1}{Corollary~\ref*{#1}}}
\newrefformat{cha}{\savehyperref{#1}{Chapter~\ref*{#1}}}
\newrefformat{sec}{\savehyperref{#1}{Section~\ref*{#1}}}
\newrefformat{app}{\savehyperref{#1}{Appendix~\ref*{#1}}}
\newrefformat{tab}{\savehyperref{#1}{Table~\ref*{#1}}}
\newrefformat{fig}{\savehyperref{#1}{Figure~\ref*{#1}}}
\newrefformat{hyp}{\savehyperref{#1}{Hypothesis~\ref*{#1}}}
\newrefformat{alg}{\savehyperref{#1}{Algorithm~\ref*{#1}}}
\newrefformat{rem}{\savehyperref{#1}{Remark~\ref*{#1}}}
\newrefformat{item}{\savehyperref{#1}{Item~\ref*{#1}}}
\newrefformat{step}{\savehyperref{#1}{step~\ref*{#1}}}
\newrefformat{conj}{\savehyperref{#1}{Conjecture~\ref*{#1}}}
\newrefformat{fact}{\savehyperref{#1}{Fact~\ref*{#1}}}
\newrefformat{prop}{\savehyperref{#1}{Proposition~\ref*{#1}}}
\newrefformat{prob}{\savehyperref{#1}{Problem~\ref*{#1}}}
\newrefformat{claim}{\savehyperref{#1}{Claim~\ref*{#1}}}
\newrefformat{relax}{\savehyperref{#1}{Relaxation~\ref*{#1}}}
\newrefformat{red}{\savehyperref{#1}{Reduction~\ref*{#1}}}
\newrefformat{part}{\savehyperref{#1}{Part~\ref*{#1}}}
\newrefformat{ex}{\savehyperref{#1}{Exercise~\ref*{#1}}}
\newrefformat{property}{\savehyperref{#1}{Property~\ref*{#1}}}

\newcommand{\Sref}[1]{\hyperref[#1]{\S\ref*{#1}}}

\usepackage{nicefrac}

\ifnum\fastmode=0
\ifnum\usemicrotype=1
\usepackage{microtype}
\fi
\fi

\ifnum\showauthornotes=2
\newcommand{\mynotes}[1]{{\sffamily\small\color{teal}{#1}}\medskip}
\newcommand{\Authornote}[2]{{\sffamily\small\color{blue}{[#1: #2]}}\medskip}
\newcommand{\Authornotecolored}[3]{{\sffamily\small\color{#1}{[#2: #3]}}}
\newcommand{\Authorcomment}[2]{{\sffamily\small\color{gray}{[#1: #2]}}}
\newcommand{\Authorstartcomment}[1]{\sffamily\small\color{gray}[#1: }

\newcommand{\Authorfnote}[2]{\footnote{\color{red}{#1: #2}}}
\newcommand{\Authorfixme}[1]{\Authornote{#1}{\textbf{??}}}
\newcommand{\Authormarginmark}[1]{\marginpar{\textcolor{red}{\fbox{\Large #1:!}}}}
\newcommand{\myexplain}[1]{{\sffamily\small\color{red}{\noindent [Explanation:\medskip\newline \begin{quote}#1\hfill]\end{quote}}}\medskip}
\newcommand{\explain}[1]{{\sffamily\small\color{red}{#1}}\medskip}

\else

\newcommand{\mynotes}[1]{}
\newcommand{\Authornote}[2]{}
\newcommand{\Authornotecolored}[3]{}
\newcommand{\Authorcomment}[2]{}
\newcommand{\Authorstartcomment}[1]{}

\newcommand{\Authorfnote}[2]{}
\newcommand{\Authorfixme}[1]{}
\newcommand{\Authormarginmark}[1]{}
\newcommand{\myexplain}[1]{}
\newcommand{\explain}[1]{}

\ifnum\showauthornotes=1
\renewcommand{\myexplain}[1]{{\sffamily\small\color{red}{\noindent \begin{quote}{\bf Explanation:} \medskip\newline #1\end{quote}}}\medskip}
\fi

\fi

\ifnum\showfixme=0

\fi

\usepackage{boxedminipage}

\newcommand{\textparen}[1]{\text{(#1)}}

\ifx\because\undefined
\newcommand{\because}[1]{\textparen{because #1}}
\else
\renewcommand{\because}[1]{\textparen{because #1}}
\fi

\newcommand{\defeq}{\stackrel{\mathrm{def}}=}

\newcommand\bdot\bullet

\ifx\mathds\undefined %
\newcommand{\Ind}{\mathbb I}
\else
\newcommand{\Ind}{\mathds 1}
\fi

\newcommand{\Z}{\mathbb Z}

\renewcommand{\leq}{\leqslant}

\renewcommand{\geq}{\geqslant}

\ifnum\showdraftbox=1

\else

\fi

\let\epsilon=\varepsilon

\numberwithin{equation}{section}

\newcommand\MYcurrentlabel{xxx}

\newcommand{\MYstore}[2]{%
  \global\expandafter \def \csname MYMEMORY #1 \endcsname{#2}%
}

\newcommand{\MYload}[1]{%
  \csname MYMEMORY #1 \endcsname%
}

\newcommand{\MYnewlabel}[1]{%
  \renewcommand\MYcurrentlabel{#1}%
  \MYoldlabel{#1}%
}

\newcommand{\MYdummylabel}[1]{}

\newcommand{\torestate}[1]{%
  \let\MYoldlabel\label%
  \let\label\MYnewlabel%
  #1%
  \MYstore{\MYcurrentlabel}{#1}%
  \let\label\MYoldlabel%
}

\newcommand{\restatetheorem}[1]{%
  \let\MYoldlabel\label
  \let\label\MYdummylabel
  \begin{theorem*}[Restatement of \prettyref{#1}]
    \MYload{#1}
  \end{theorem*}
  \let\label\MYoldlabel
}

\newcommand{\restatelemma}[1]{%
  \let\MYoldlabel\label
  \let\label\MYdummylabel
  \begin{lemma*}[Restatement of \prettyref{#1}]
    \MYload{#1}
  \end{lemma*}
  \let\label\MYoldlabel
}

\newcommand{\restateprop}[1]{%
  \let\MYoldlabel\label
  \let\label\MYdummylabel
  \begin{proposition*}[Restatement of \prettyref{#1}]
    \MYload{#1}
  \end{proposition*}
  \let\label\MYoldlabel
}

\newcommand{\restatefact}[1]{%
  \let\MYoldlabel\label
  \let\label\MYdummylabel
  \begin{fact*}[Restatement of \prettyref{#1}]
    \MYload{#1}
  \end{fact*}
  \let\label\MYoldlabel
}

\newcommand{\restate}[1]{%
  \let\MYoldlabel\label
  \let\label\MYdummylabel
  \MYload{#1}
  \let\label\MYoldlabel
}

\newcommand{\addreferencesection}{
  \phantomsection
\ifnum\stocmode=0
  \addcontentsline{toc}{section}{References}
\else
  \addcontentsline{toc}{section}{References \hspace*{1in} --------- End of extended abstract ---------}
\fi

}

\newcommand{\e}{\epsilon}
\newcommand{\eps}{\epsilon}

\let\origparagraph\paragraph
\renewcommand{\paragraph}[1]{\vspace*{-2pt}\origparagraph{#1.}}

\allowdisplaybreaks

\usepackage{paralist}

\usepackage{comment}

\let\pref=\prettyref

\newcommand{\vertiii}[1]{{\left\vert\kern-0.25ex\left\vert\kern-0.25ex\left\vert #1 
          \right\vert\kern-0.25ex\right\vert\kern-0.25ex\right\vert}}

\renewcommand{\Ind}{\vvmathbb{1}}
\ifnum\jamesmode=0
\fi
\let\1\Ind

\renewcommand{\Ind}{\mathbf{1}}

\renewcommand{\1}{\mathbf{1}}

\renewcommand{\Z}{\vvmathbb{Z}}

\begin{document}

\title{On expanders from the action of $GL(2,\vvmathbb{Z})$}
\author{James R. Lee\thanks{Department of Computer Science \& Engineering, University of Washington.  Partially supported by NSF grants CCF-1217256 and CCF-0915251.} }
\date{}

\maketitle

\begin{abstract}
Consider the undirected graph $G_n=(V_n, E_n)$ where $V_n = (\Z/n\Z)^2$ and
$E_n$ contains an edge from $(x,y)$ to $(x+1,y)$, $(x,y+1)$, $(x+y,y)$, and $(x,y+x)$
for every $(x,y) \in V_n$.  Gabber and Galil, following Margulis, gave an elementary
proof that $\{G_n\}$ forms an expander family.  In this expository note, we present a somewhat simpler proof
of this fact, and demonstrate its utility by
isolating a key property of the linear transformations $(x,y) \mapsto (x+y,x), (x,y+x)$
that yields expansion.

As an example, take any invertible, integral matrix $S \in GL_2(\Z)$ and let
$G^{S}_n = (V_n, E^{S}_n)$ where $E^{S}_n$ contains,
for every $(x,y) \in V_n$, an edge from $(x,y)$ to $(x+1,y)$, $(x,y+1)$, $S(x,y)$, and $S^{\top}(x,y)$,
and $S^{\top}$ denotes the transpose of $S$.
Then $\{G_n^{S}\}$ forms an expander family if and only if the infinite
graph $$G^{S} = \left(\Z^2 \setminus \{0\}, \left\{ \vphantom{\bigoplus}\{z,S z\}, \{z, S^{\top} z\} : z \in \Z^2 \setminus \{0\}\right\}\right)$$
has positive Cheeger constant.

This latter property turns out to be elementary to analyze:
For any $S = \left(\begin{smallmatrix} a & b \\ c & d \end{smallmatrix}\right)\in GL_2(\Z)$, the graph $G^S$ has positive Cheeger constant if and
only if $(a+d)(b-c) \neq 0$.
The case $S = \left(\begin{smallmatrix} 1 & 1 \\ 0 & 1 \end{smallmatrix}\right)$
recovers the Margulis-Gabber-Galil graphs.  We also present some other generalizations.
\end{abstract}

\section{Introduction}

Expander graphs have played a fundamental role in many areas of mathematics and computer science;
we refer to the monograph \cite{HLW06}.  Margulis \cite{Margulis73} discovered
 the first explicit construction of expanders.  Based on his work, Gabber and Galil \cite{GG81}
later presented an elementary construction and analysis.  The Gabber-Galil graphs
still provide the simplest, most succinct description of expanders to date.

Consider the undirected graph $G_n=(V_n,E_n)$ with vertex set
$V_n = (\Z/n\Z)^2$ and edge set $E_n$
which contains, for every $(x,y) \in V_n$, an edge to each of
$(x\pm 1,y), (x,y\pm 1), (x\pm y,y), (x,y\pm x)$.
Then $\{G_n : n \geq 2\}$ forms a family of expander graphs with vertex degree at most $8$.
Jimbo and Maruoka \cite{JM87}, using discrete Fourier analysis, presented another proof that the Gabber-Galil graphs
are expanders.  Both these analyses contain at least
one non-trivial and arguably opaque technical analytic step.
For instance, the survey \cite{HLW06} gives an elementary proof along the lines of \cite{JM87}
but still refers to the argument as ``subtle and mysterious.''

We present a somewhat simpler proof, or at least one
whose pieces are each well-motivated.  The ``technical step'' is replaced by an application
of the discrete Cheeger inequality and a very simple combinatorial lemma inspired by a paper of
Linial and London \cite{LL06} (cf. \pref{lem:dvir}).
Moreover, the basic approach allows us to analyze a variety of similar families.

Given any two invertible, integral matrices $S,T \in GL_2(\Z)$, one can consider
the family of graphs $G_n^{S,T} = (V_n, E_n^{S,T})$, where $E_n^{S,T}$
contains edges from every $(x,y) \in V_n$ to each of \[(x \pm 1,y), (x,y\pm 1), S(x,y), S^{-1}(x,y), T(x,y), T^{-1}(x,y)\,.\]
The Gabber-Galil graphs correspond to the choice $S = \left(\begin{smallmatrix} 1 & 1 \\ 0 & 1 \end{smallmatrix}\right)$ and $T = \left(\begin{smallmatrix} 1 & 0 \\ 1 & 1 \end{smallmatrix}\right)$.

Consider also the countably infinite graph
$G^{S,T}$ with vertex set $\Z^2 \setminus \{0\}$
and edges $$E^{S,T} \defeq \left\{ \{z,S z\}, \{z, T z\} : z \in \Z^2 \setminus \{0\} \right\}\,.$$
In \pref{sec:general}, we prove the following relationship.

\begin{theorem}\label{thm:intro1}
For any $S,T \in GL_2(\Z)$, if $G^{S^{\top},T^{\top}}$ has positive Cheeger constant, then $\{G_n^{S,T}\}$
is a family of expander graphs.
\end{theorem}

An infinite graph $G=(V,E)$ with uniformly
bounded degrees has positive
Cheeger constant if there is a number $\e > 0$ such that every finite subset $U \subseteq V$
has at least $\e |U|$ edges with exactly one endpoint in $U$.
While \pref{thm:intro1} may not seem particularly powerful,
it turns out that in many interesting cases, proving a non-trivial lower bound on
the Cheeger constant of $G^{S,T}$ is elementary.
For the Gabber-Galil graphs, the argument is especially simple;
see \pref{lem:dvir}.

One can generalize the Gabber-Galil graphs in a few different ways.
As a prototypical example, consider the family $\{G_n^{S,S^{\top}}\}$
for any $S \in GL_2(\Z)$.  In \pref{sec:expansion},
we give the following characterization.

\begin{theorem}
For any $S =\left(\begin{smallmatrix} a & b \\ c & d \end{smallmatrix}\right) \in GL_2(\Z)$, it holds that $\{G_n^{S,S^{\top}}\}$ is an
expander family if and only if $(a+d)(b-c) \neq 0$.
\end{theorem}

For instance, the preceding theorem implies that if $S$ has order $4$ then $\{G_n^{S,S^{\top}}\}$ is
not a family of expander graphs, but if $S$ has order $6$  and $S \neq S^{\top}$ then the graphs
are expanders.

Earlier, Cai \cite{Cai03} considered a different generalization.  Let $R = \left(\begin{smallmatrix} 0 & 1 \\ 1 & 0 \end{smallmatrix}\right)$ be the reflection across the line $y=x$.
The Gabber-Galil graphs can also be seen as $G_n^{S,T}$ where $S = \left(\begin{smallmatrix} 1 & 1 \\ 0 & 1 \end{smallmatrix}\right)$ and $T = RSR$.
In \pref{sec:cai}, we give the following characterization.

\begin{theorem}
For any $S =\left(\begin{smallmatrix} a & b \\ c & d \end{smallmatrix}\right) \in GL_2(\Z)$, it holds that $\{G_n^{S,RSR}\}$ is an
expander family if and only if $(a+d)(b+c) \neq 0$.
\end{theorem}

Cai \cite{Cai03} considers the situation $\det(S)=1$ and $|a+d| \geq 2, |b+c| \geq 2$.
However, his work does not prove that $\{G_n^{S,RSR}\}$ are expanders.
In fact, the graphs he associates to a matrix $S$ are somewhat complicated
and need to refer to the action of $S$ on the torus.  Moreover, they
do not have uniformly bounded degree; the degree of his graphs grow linearly in $\|S\|_1$
(the sum of the magnitudes of the entries of $S$).
The maximum degree of our graphs is clearly bounded by 8.
Interestingly, Cai states that $\{G_n^{S,S^{\top}}\}$ is
a more natural generalization, but the main technical tool
of the Gabber-Galil style analysis (see \pref{thm:cai})
does not work for these graphs.

\section{The Margulis-Gabber-Galil graphs}
\label{sec:GG}

Consider an undirected graph $G=(V,E)$ with an at most countable vertex set.
For $A,B \subseteq V$, we use $E(A,B)$ to denote the
set of edges with one endpoint in $A$ and one in $B$.
We write $E(A) = E(A,\bar A)$ where $\bar A$ denotes the complement of $A$ in $V$.
We define the expansion of a subset $U \subseteq V$ by
$$
h_G(U) \defeq \frac{|E(U)|}{|U|}\,.
$$
For $G$ finite, we set $h(G) \defeq \min_{|U| \leq \frac12 |V|} h_G(U)$.
If $G$ is infinite, we put $h(G) \defeq \min_{U \subseteq V : |U| < \infty} h_G(U)$.
In both the finite and infinite case, we refer to $h(G)$ as the {\em Cheeger constant of $G$.}

We also have the Rayleigh quotient of a function $f : V \to \vvmathbb{C}$ given by
$$
\mathcal R_G(f) \defeq \frac{\sum_{\{u,v\}\in E} |f(u)-f(v)|^2}{\sum_{u \in V} |f(u)|^2}\,,
$$
and for finite $G$, we put $\lambda_2(G) \defeq \min \{ \mathcal R_G(f) : \sum_{u \in V} f(u) = 0 \}$.
This is the smallest non-zero eigenvalue of the combinatorial Laplacian (see, e.g., the book \cite{Chung97}).
An infinite family of finite graphs $\{G_n\}$ with uniformly bounded degrees
is called an {\em expander family} if $\lambda_2(G_n) \geq c > 0$ for some $c > 0$.
We will assume familiarity with the following discrete Cheeger inequality.

\begin{lemma}\label{lem:cheeger}
For any countable graph $G=(V,E)$ with maximum degree $\Delta$ and any function $f : V \to \vvmathbb{C}$ with
$\sum_{v \in V} |f(v)|^2 < \infty$, there exists a
finite subset $U \subseteq \{ v \in V : f(v) \neq 0\}$ such that
$$
h_G(U) \leq \sqrt{2 \Delta \mathcal R_G(f)}\,.
$$
\end{lemma}

\begin{proof}
Let $U_t = \{ v \in V : |f(v)|^2 \geq t \}$.
Observe that for each $t > 0$, one has $U_t \subseteq \{ v \in V : f(v) \neq 0\}$ and $U_t$
is finite since $\sum_{v \in V} |f(v)|^2$ is finite.
Now we have:
\begin{eqnarray*}
\int_0^{\infty} |E(U_t, \bar U_t)|\,dt &=& \sum_{\{u,v\} \in E} \left||f(u)|^2 - |f(v)|^2\right| \\
&= &
\sum_{\{u,v\} \in E} (|f(u)| + |f(v)|) (|f(u)|-|f(v)|) \\
&\leq &
\sqrt{\sum_{\{u,v\} \in E} (|f(u)|+|f(v)|)^2} \sqrt{\sum_{\{u,v\} \in E} |f(u)-f(v)|^2} \\
&\leq &
\sqrt{2 \Delta \sum_{u \in V} |f(u)|^2}\sqrt{\sum_{\{u,v\} \in E} |f(u)-f(v)|^2}.
\end{eqnarray*}
On the other hand,
$\int_0^{\infty} |U_t|\,dt = \sum_{u \in V} |f(u)|^2,$
thus
$$
\int_0^{\infty} |E(U_t, \bar U_t)|\,dt \leq \sqrt{2\Delta\mathcal R_G(f)} \int_0^{\infty} |U_t|\,dt\,,
$$
implying there exists a $t > 0$ such that $h_G(U_t) \leq \sqrt{2\Delta\mathcal R_G(f)}$.
\end{proof}

\medskip
\noindent
{\bf An initial expanding object.}
We will start with an initial ``expanding object,'' and then try
to construct a family of graphs out of it.  First, consider the infinite graph
$\mathcal G=(\Z^2, E)$ whose edges are given by two maps $S, T : \vvmathbb{R}^2 \to \vvmathbb{R}^2$ defined
by $S(x,y) = (x,x+y)$ and $T(x,y) = (x+y,y)$.  Each vertex $z \in \Z^2$
is connected to $S(z), S^{-1}(z), T(z), T^{-1}(z)$.  So every vertex has degree
at most four.  Clearly $(0,0)$ is not adjacent to anything.
Using an argument from \cite{LL06}, we will show that this graph is an expander in the following sense.

\begin{lemma}\label{lem:dvir}
For any finite subset $A \subseteq \Z^2 \setminus \{0\}$, we have
$
|E(A, \bar A)| \geq |A|
$.
\end{lemma}

\begin{proof}
Define $Q_1 = \{ (x,y) \in \Z^2 : x > 0, y \geq 0 \}$.  This is the first quadrant, without the $y$-axis and the origin.
Define $Q_2, Q_3, Q_4$ similarly by rotating $Q_1$ by $90$, $180$, and $270$ degrees, respectively,
and note that we have a partition $\Z^2 \setminus \{0\} = Q_1 \cup Q_2 \cup Q_3 \cup Q_4$.

Let $A_i = A \cap Q_i$.  We will show that $|E(A_1, \bar A \cap Q_1)| \geq |A_1|$.  Since our graph is invariant
under rotations of the plane by $90^{\circ}$, this will imply our goal:
$$
|E(A, \bar A)| \geq \sum_{i=1}^4 |E(A_i, \bar A \cap Q_i)| \geq \sum_{i=1}^4 |A_i| = |A|\,.
$$

It is immediate that $S(A_1), T(A_1) \subseteq Q_1$.  Furthermore, we have $S(A_1) \cap T(A_1) = \emptyset$ because
$S$ maps points in $Q_1$ above (or onto) the line $y=x$ and $T$ maps points of $Q_1$ below the line $y=x$.
Furthermore, $S$ and $T$ are bijections, thus $|S(A_1) + T(A_1)| = |S(A_1)| + |T(A_1)| = 2|A_1|\,.$
In particular, this yields $|E(A_1, \bar A \cap Q_1)| \geq |A_1|$, as desired.
\end{proof}

Of course, $\mathcal G$ is not a finite graph, so for a number $n \geq 2$, we define the
graph $G_n=(V_n,E_n)$ with vertex set $V_n = (\Z/n \Z)^2$.
There are four types of edges in $E_n$:  A vertex $(x,y)$ is connected to the
vertices $$\{(x,y\pm 1), (x\pm 1, y), (x,x\pm y), (x\pm y,y)\}\,,$$
where arithmetic is taken modulo $n$.
This yields a graph of degree at most 8.
We now state the main result of this section.

\begin{theorem}
\label{thm:main}
There is a constant $c > 0$ such that for every $n \geq 2$,
$$
\lambda_2(G_n) \geq c\,.
$$
In other words, $\{G_n\}$ forms an expander family.
\end{theorem}

\medskip
\noindent
{\bf Passing to the continuous torus.}
Our results for
$\Z^2$ do not seem immediately useful for analyzing these finite graphs.
We will first pass from the discrete graphs $\{G_n\}$
to the continuous torus.  This is a reassuring step, as it means
our analysis is not going to rely on number theoretic considerations of the modulus $n$.

Let $\vvmathbb{T}^2 = \vvmathbb{R}^2/\Z^2$ be the 2-dimensional torus equipped
with the Lebesgue measure and consider the complex Hilbert space
$$
L^2(\vvmathbb{T}^2) = \left\{ f : \vvmathbb{T}^2 \to \vvmathbb{C} : \int_{\vvmathbb{T}^2} |f|^2 < \infty \right\}\,.
$$
equipped with the inner product $\langle f,g\rangle_{L^2} = \int_{\vvmathbb{T}^2} f \bar g$.

\medskip

We also define a related value
\begin{equation}
\label{eq:L2T2}
\lambda_2(\vvmathbb{T}_{S,T}^2) \defeq
\min_{f \in L^2(\vvmathbb{T}^2)} \left\{\frac{\|f-f\circ S\|_{L^2}^2 + \|f-f \circ T\|_{L^2}^2}{\|f\|_{L^2}^2} : \int_{\vvmathbb{T}^2} f = 0\right\}\,.
\end{equation}

\begin{lemma}\label{lem:torus}
There is some $\varepsilon > 0$ such that
for any $n \geq 2$, we have $\lambda_2(G_n) \geq \varepsilon \lambda_2(\vvmathbb{T}_{S,T}^2)$\,.
\end{lemma}

\begin{proof}
Suppose we are given some map $f : V_n \to \vvmathbb{C}$ such that $\sum_{u \in V_n} f(u)=0$.
We define its continuous extension $\tilde f : \vvmathbb{T}^2 \to \vvmathbb{C}$ as
follows.  There is a natural embedding of $V_n$ into $[0,1]^2$ which
we represent as follows:  Given a point $w=(x/n,y/n) \in [0,1]^2$,
with $x,y \in \{0,1,\ldots,n\}$, we write $[w]$ for the corresponding element of $V_n$.

Every point $z \in [0,1)^2$ sits inside a grid square with four
corners $u_1,u_2,u_3,u_4$ such that $[u_1],[u_2],[u_3],[u_4] \in V_n$.
We call such a square (thought of as a subset of $\vvmathbb{T}^2$) a {\em canonical square.}
Define $\tilde f(z)$ as the average
\begin{equation}\label{eq:extension}
\tilde f(z) = \frac{\sum_{i=1}^4 (\frac1{n}-\|u_i-z\|_{\infty}) f([u_i])}
                   {\sum_{i=1}^4 (\frac1{n}-\|u_i-z\|_{\infty})}\,.
\end{equation}
Observe that this is well-defined; e.g., if $z$ lies on the segment between $u_1$ and $u_2$
then the coefficients of $f([u_3])$ and $f([u_4])$ are zero.
By symmetry, it follows immediately that $\int_{\vvmathbb{T}^2} \tilde f = 0$.

It is also easy to verify that $\|\tilde f\|^2_{L^2} \geq \frac{c}{n^2} \sum_{v \in V} f(v)^2$
for some $c > 0$.  For any square with corners $\{u_1,u_2,u_3,u_4\}$, let $i \in \{1,2,3,4\}$
be such that $f([u_i])^2$ is maximal and let $B$ denote an $\ell_{\infty}$ ball of radius $\frac1{8n}$
around $u_i$.  Then $\int_B |\tilde f|^2 \geq \frac{c}{n^2} \sum_{i=1}^4 f([u_i])^2$ for some universal constant $c > 0$.
Summing over all the squares yields the claim.

So to finish the proof,
we are left to argue that
\begin{equation}\label{eq:grad}
\|\tilde f-\tilde f\circ S\|_{L^2}^2 + \|\tilde f-\tilde f\circ S\|_{L^2}^2 \leq \frac{c}{n^2} \sum_{\{u,v\} \in E_n} (f(u)-f(v))^2
\end{equation}
for some $c > 0$.  Consider any point $z \in \vvmathbb{T}^2$ contained in a square $\square_1$ and suppose $S(z)$ is in $\square_2$.
Note that $\square_1=\square_2$ is a possibility.
Let $\mathcal C$ be the set of (at most) eight vertices of $V_n$ that comprise the corners of $\square_1$ and $\square_2$.
Then any pair of vertices in $\mathcal C$ can reach each other using a path of length at most five in $G_n$.
This is the only place where we need to use
the fact that edges of the form $(x,y) \leftrightarrow (x,y\pm 1)$ and $(x,y) \leftrightarrow (x\pm 1,y)$
are present in $G_n$.
On the other hand, we clearly have
$$
|\tilde f(z) - \tilde f(S(z))|^2 \leq \max_{u,v \in \mathcal C} |f(u)-f(v)|^2\,,
$$
since $\tilde f(z)$ is a convex combination of the $f$-values at the corners of $\square_1$
and $\tilde f(S(z))$ is a convex combination of the $f$-values at the corners of $\square_2$.

Now consider a canonical square $\square \subseteq \vvmathbb{T}^2$,
which has measure $1/n^2$.
Let $E(\square)$ to denote the set of edges in $G_n$ that occur on some path of length at most 5
emanating from the corners of $\square$.  Then the preceding argument yields
$$
\int_{\square} |\tilde f(z) - \tilde f(S(z))|^2 dz \leq \frac{1}{n^2} \max_{\{u,v\} \in E(\square)} |f(u)-f(v)|^2 \leq O\left(\frac{1}{n^2}\right) \sum_{\{u,v\} \in E(\square)} |f(u)-f(v)|^2\,,
$$
using the fact that $|E(\square)| = O(1)$ because $G_n$ has degree at most $8$.  Summing the preceding inequality over all canonical squares yields
$$
\int_{\vvmathbb{T}^2}  |\tilde f(z) - \tilde f(S(z))|^2 dz \leq O\left(\frac{1}{n^2}\right) \sum_{\{u,v\} \in E} |f(u)-f(v)|^2\,,
$$
since every edge occurs in some set $E(\square)$ at most $O(1)$ times.  An identical argument holds for $T$, yielding
\eqref{eq:grad}.
\end{proof}

\medskip
\noindent
{\bf Using the Fourier transform to unwrap the torus.}
Our final goal is to show that $\lambda_2(\vvmathbb{T}_{S,T}^2) > 0$.
Our approach is based on the fact that $S$ and $T$, being shift operators,
will act rather nicely on the Fourier basis.

We recall that if $m,n \in \vvmathbb{N}$ and we define $\chi_{m,n} \in L^2(\vvmathbb{T}^2)$ by
$\chi_{m,n}(x,y) = \exp(2\pi i(mx+ny))$, then $\{\chi_{m,n} : m,n \in \Z\}$ forms an orthonormal
Hilbert basis for $L^2(\vvmathbb{T}^2)$.  In particular, every $f \in L^2(\vvmathbb{T}^2)$
can be written as
\begin{equation}\label{eq:as}
f = \sum_{m,n \in \Z} \hat f(m,n) \chi_{m,n}\,,
\end{equation}
where $\hat f(m,n) = \langle f, \chi_{m,n}\rangle_{L^2}$
and convergence in \eqref{eq:as} is in the $L^2(\vvmathbb{T}^2)$ norm (see, for instance, \cite[\S I.5]{Katznelson04}).
Putting $\ell^2(\Z^2) = \{f : \Z^2 \to \vvmathbb{C} : \sum_{z \in \Z^2} |f(z)|^2 < \infty \}$,
the Fourier transform is the linear isometry $f \mapsto \hat f$
from $L^2(\vvmathbb{T}^2)$ to $\ell^2(\Z^2)$.

For any $m,n \in \Z$, we have
$$
\chi_{m,n} \circ S = \chi_{m,n+m} \quad\textrm{ and }\quad \chi_{m,n} \circ T= \chi_{m+n,n}.
$$
Thus for any $f \in L^2(\vvmathbb{T}^2)$,
we have
\begin{eqnarray*}
\widehat{f \circ S} &=&  \sum_{m,n} \hat f(m,n) \chi_{m,n+m} = \sum_{m,n} \hat f(m,n-m) \chi_{m,n} = \hat f \circ T^{-1} \\
\widehat{f \circ T} &=& \sum_{m,n} \hat f(m,n) \chi_{m+n,n} = \sum_{m,n} \hat f(m-n,n) \chi_{m,n} = \hat f \circ S^{-1}\,.
\end{eqnarray*}
The final thing to note is that $\hat f(0,0) = \langle f, \chi_{0,0} \rangle = \int_{\vvmathbb{T}^2} f$.
So now if we simply apply the Fourier transform (a linear isometry) to the expression
in \eqref{eq:L2T2}, we arrive at
\begin{eqnarray*}
\lambda_2(\vvmathbb{T}_{S,T}^2) &=&
\min_{f \in L^2(\vvmathbb{T}^2)} \left\{\frac{\|\hat f-\widehat{f\circ S}\|_{\ell^2}^2 + \|\hat f-\widehat{f \circ T}\|_{\ell^2}^2}{\|\hat f\|_{\ell^2}^2} : \int_{\vvmathbb{T}^2} f = 0\right\} \\
&=&
\min_{\hat f \in \ell^2(\Z^2)}
\left\{\frac{\sum_{z \in \Z^2} |\hat f(z)- \hat f(T^{-1}(z))|^2 + |\hat f(z)-\hat f(S^{-1}(z))|^2}{\|\hat f\|_{L^2}^2} : \hat f(0,0)=0\right\}\,.
\end{eqnarray*}

In other words, $$\lambda_2(\vvmathbb{T}_{S,T}^2) = \min_{\hat f \in \ell^2(\Z^2)} \{ \mathcal R_{\mathcal G}(\hat f) : \hat f(0,0) = 0\}\,,$$
where $\mathcal G$ is our initial graph defined on $\Z^2$.
Applying the discrete Cheeger inequality (\pref{lem:cheeger}) with $\Delta=4$, yields
$$
\min_{\hat f : V \to \vvmathbb{C}} \{ \mathcal R_{\mathcal G}(\hat f) : \hat f(0,0) = 0 \} \geq \frac1{8} \min_{U : (0,0) \notin U} h_{\mathcal G}(U)^2 \geq \frac1{8}\,,
$$
where the final inequality is exactly the content of \pref{lem:dvir}.
Thus by \pref{lem:torus} for some $\e > 0$ and every $n \geq 2$, we have $\lambda_2(G_n) \geq \e \lambda_2(\vvmathbb{T}_{S,T}^2) \geq \frac{\e}8$.
This completes the proof of \pref{thm:main}.

\section{The general correspondence}
\label{sec:general}

We now perform the steps of the preceding section
is somewhat greater generality.
Consider $S,T \in GL_2(\Z)$.
We will write $\hat G^{S,T}$ to denote $G^{S^{\top},T^{\top}}$.
The main result of this section is a connection between the expansion of $\{G_n^{S,T}\}$ and $\hat G^{S,T}$.

\begin{theorem}\label{thm:maingen}
For every $S,T \in GL_2(\Z)$, if $h(\hat G^{S,T}) > 0$, then
$\{G_n^{S,T}\}$ forms an expander family.
\end{theorem}

Define the quantity
\[
\lambda_2(\vvmathbb{T}_{S,T}^2) \defeq \min_{f \in L^2(\vvmathbb{T}^2)} \left\{ \frac{\|f-f\circ S\|_{L^2}^2 + \|f-f\circ T\|_{L^2}2}{\|f\|_{L^2}}
: \int_{\vvmathbb{T}^2} f = 0 \right\}\,.
\]

The following result requires a bit more delicacy than \pref{lem:torus}.

\begin{lemma}\label{lem:torusgen}
There is an $\e > 0$ such that for every $S,T \in GL_2(\Z)$ and $n \geq 2$, we have
$$
\lambda_2(G_n^{S,T}) \geq \frac{\e}{\|S\|_1^2 + \|T\|_1^2}\, \lambda_2(\vvmathbb{T}_{S,T}^2)\,.
$$
\end{lemma}

\begin{proof}
We will use the notion of canonical squares from \pref{lem:torus}.
Suppose we have a map $f : V_n \to \vvmathbb{C}$ satisfying $\sum_{u \in V_n} f(u)=0$.
Define the extension $\tilde f : \vvmathbb{T}^2 \to \vvmathbb{C}$ as in \eqref{eq:extension}.  The fact that $\int_{\vvmathbb{T}^2} \tilde f  = 0$
and $\int_{\vvmathbb{T}^2} |\tilde f|^2 \geq \frac{c}{n^2} \sum_{u \in V_n} |f(u)|^2$ for some absolute constant $c > 0$ is proved in \pref{lem:torus}.
We are thus left to prove that for some $c > 0$,
\begin{equation}\label{eq:gradient}
\|\tilde f-\tilde f \circ S\|_{L^2}^2 + \|\tilde f-\tilde f\circ T\|_{L^2}^2 \leq c\frac{\|S\|^2_1+\|T\|^2_1}{n^2} \sum_{\{u,v\} \in E^{S,T}_n} |f(u)-f(v)|^2\,.
\end{equation}

To this end, suppose $S = \left(\begin{smallmatrix} a & b \\ c & d \end{smallmatrix}\right)$
and consider a point $z \in \square_1$ and $S z \in \square_2$,
where $\square_1$ and $\square_2$ are canonical squares whose corners are vertices from $V_n$ (it is possible that $\square_1=\square_2$).  Since $\tilde f(z)$ is a convex combination
of the values of $f$ at the corners of $\square_1$ and similarly for $\tilde f(S(z))$ and $\square_2$, we have
\begin{equation}\label{eq:max}
|\tilde f(z) - \tilde f(Sz)|^2 \leq \max_{u,v \in \mathcal C} |f(u)-f(v)|^2\,,
\end{equation}
where $\mathcal C$ contains the (at most) eight corners of $\square_1$ and $\square_2$.

Unlike in \pref{lem:torus}, the members of $\mathcal C$ can no longer be connected by paths of length $O(1)$ in $G_n^S$.
However, it is elementary to see that they can be connected by paths of length at most $\|S\|_1+1=|a|+|b|+|c|+|d|+1$.
We simply need to choose the paths in a consistent way in order to conclude that \eqref{eq:gradient} holds.
This will be a bit technical, but the underlying idea is very simple.

We will now specify canonical paths between the members of $\mathcal C$.
Let us write $E'_n \subseteq E^{S,T}_n$ for the set of edges connecting $(x,y)$ to $(x\pm 1,y)$ or $(x,y\pm 1)$.
Call an edge of $E'_n$ {\em horizontal} if it changes the $x$ coordinate and {\em vertical} if it changes the $y$ coordinate.

Let $(x,y) \in [0,1)^2$ denote the lower-left corner of $\square_1$ and let $(x',y') \in [0,1)^2$ denote the lower-left corner of $\square_2$.
We may assume that $z=(x+\alpha,y+\beta)$ for some $\alpha,\beta \in (0,1/n)$,
and \[Sz = S(x,y) + S(\alpha,\beta) = S(x,y) + (a\alpha + b\beta, c\alpha + d\beta)\,.\]
We specify a path from $(x,y)$ to $(x',y')$.
Our path $P_z$ in $G_n^S$ will first follow the edge $\{(x,y), S(x,y)\}$
then move along edges of $E'_n$ in the $x$ direction for $\lfloor a\alpha + b \beta\rfloor$ steps,
then move along edges of $E'_n$ in the $y$ direction for $\lfloor c\alpha + d\beta \rfloor$ steps.
This will arrive at some corner of $\square_2$ (e.g., the lower-left corner if all the entries of $S$ are positive).
Our path then moves to $(x',y')$ using at most two additional edges of $\square_2$.
For any other pair $u,v \in \mathcal C$:  If they are in the same square, move along the edges of the square
in some canonical way using a path of length at most two.  Otherwise, if $u$ is a corner of $\square_1$
and $v$ is a corner of $\square_2$, first from $u$ to $(x,y)$ along edges of $\square_1$, then to $(x',y')$ using $P_z$,
then from $(x',y')$ to $v$ using edges of $\square_2$.  Let $P^z_{uv}$ denote the specified path between $u,v \in \mathcal C$.
Note that the length of $P^z_{uv}$ is $O(\|S\|_1)$.

The main points of this construction are as follows.  First, for every pair of horizontal (respectively, vertical) edges $e,e' \in E'_n$, we have
\begin{equation}\label{eq:equitable}
\int_{\vvmathbb{T}^2} \1_{\{e \in P_z\}} dz = \int_{\vvmathbb{T}^2} \1_{\{e' \in P_z\}} dz\,.
\end{equation}
The second is that, combining \eqref{eq:max} with Cauchy-Schwarz yields
\begin{equation}\label{eq:cs}
|\tilde f(z) - \tilde f(S(z))|^2 \leq O(\|S\|_1) \sum_{u,v \in \mathcal C} \sum_{\{r,s\} \in P^z_{uv}} |f(r)-f(s)|^2\,.
\end{equation}
Using the equitable property \eqref{eq:equitable} and the fact that every edge of the form $\{(x,y),S(x,y)\}$
appears on the right-hand side of \eqref{eq:cs} only when $z \in \square_1$, we can integrate \eqref{eq:cs} to yield
$$
\int_{\vvmathbb{T}^2} |\tilde f(z) - \tilde f(S(z))|^2dz \leq O\left(\frac{\|S\|^2_1}{n^2}\right) \sum_{\{u,v\} \in E^{S,T}_n} |f(u)-f(v)|^2\,.
$$
An identical analysis holds for $T$, allowing us to verify \eqref{eq:gradient}.
\end{proof}

\begin{lemma}\label{lem:fouriergen}
For any $S,T \in GL_2(\Z)$, we have
$$
\lambda_2(\vvmathbb{T}_{S,T}^2) = \min_{\hat f \in \ell^2(\Z^2)} \{ \mathcal R_{\hat G^{S,T}}(\hat f) : \hat f(0,0) = 0 \}\,.
$$
\end{lemma}

\begin{proof}
Note that if $f \in L^2(\vvmathbb{T}^2)$, then
\begin{eqnarray*}
\widehat{f \circ S} &=& \sum_{m,n} \hat f(m,n) \chi_{am+cn, bm+dn}  \\
&=& \sum_{m,n} \hat f(m,n) \chi_{S^{\top}(m,n)} \\
&=& \sum_{m,n} \hat f(S^{-\top}(m,n)) \chi_{m,n} \\
&=& \hat f \circ S^{-\top}\,.
\end{eqnarray*}
Similarly, $\widehat{f \circ T} = \hat f \circ T^{-\top}$.  Using the fact that the
Fourier transform is a linear isometry from $L^2(\vvmathbb{T}^2)$ to $\ell^2(\Z^2)$ and $\hat f(0,0) = \int_{\vvmathbb{T}^2} f$, we have
\begin{eqnarray*}
\lambda_2(\vvmathbb{T}^2_{S,T}) &=& \min_{\hat f \in \ell^2(\Z^2)} \left\{\frac{\sum_{z \in \Z^2} |\hat f(z)-\hat f(S^{-\top} z)|^2 + |\hat f(z)-\hat f(T^{-\top} z)|^2}{\sum_{z \in \Z^2} |\hat f(z)|^2} : \hat f(0,0)=0\right\} \\
&=& \min_{\hat f \in \ell^2(\Z^2)} \{ \mathcal R_{\hat G^{S,T}}(\hat f) : \hat f(0,0) = 0 \}\,,
\end{eqnarray*}
completing the proof.
\end{proof}

Combining \pref{lem:fouriergen} with the discrete Cheeger inequality (\pref{lem:cheeger}) yields the following.

\begin{corollary}
For any $S,T \in GL_2(\Z)$, $\lambda_2(\vvmathbb{T}_{S,T}^2) \geq \frac18 h(\hat G^{S,T})^2$.
\end{corollary}

Finally, combining this corollary with \pref{lem:torusgen} yields \pref{thm:maingen}.

\section{Expansion analysis}
\label{sec:expansion}

For ease of notation, we will write $G^S_n \defeq G^{S,S^{\top}}_n$ and $G^S \defeq G^{S,S^{\top}}$.

\begin{theorem}\label{thm:GS}
For any $S \in GL_2(\Z)$, it holds that $h(G^S) > 0$ if and only if
$S \neq S^{\top}$ and $\mathrm{tr}(S) \neq 0$.
\end{theorem}

Combining the preceding result with \pref{thm:maingen}, we can prove the following.

\begin{theorem}
For any $S \in GL_2(\Z)$,  it holds that $\{G_n^S\}$ is an expander family if and only if $S \neq S^{\top}$ and $\mathrm{tr}(S) \neq 0$.
\end{theorem}

\begin{proof}
Since $G^S = \hat G^{S,S^{\top}}$ and $h(G^S) > 0$ by \pref{thm:GS}, we can use \pref{thm:maingen}
to conclude that $\{G_n^S\}$ is an expander family.  On the other hand, if $S = S^{\top}$, then \pref{lem:onefin}
shows that $\{G_n^S\}$ is not an expander family.  If $\mathrm{tr}(S)=0$ then $S^4 = I = (S^{\top})^4$
and \pref{lem:fourthfin} shows that $\{G_n^S\}$ is not an expander family.
\end{proof}

To prove \pref{thm:GS}, we will first analyze the case when $\det(S)=1$ and $S$ has all non-negative entries.
This is essentially the main technical lemma of the section;
we will show that all other cases can be reduced to this one.

\begin{lemma}\label{lem:transpose}
If $S \in GL_2(\Z)$ has all non-negative entries, $\det(S)=1$, and $S \neq S^{\top}$, then
\begin{eqnarray*}
S(Q_1) \cap S^{\top}(Q_1) &=& \emptyset \\
S(Q_3) \cap S^{\top}(Q_3) &=& \emptyset \\
S^{-1}(Q_2) \cap S^{-\top}(Q_2) &=& \emptyset \\
S^{-1}(Q_4) \cap S^{-\top}(Q_4) &=& \emptyset
\end{eqnarray*}
\end{lemma}

\begin{proof}
Let $S = \left(\begin{smallmatrix} a&b \\ c &d\end{smallmatrix}\right)$ for some $a,b,c,d \geq 0$ and let $T = S^{\top}$.
Since $\det(S)=1$,
we can write:
\begin{equation}\label{eq:inverses}
S^{-1} = \left(\begin{smallmatrix} d&-b \\ -c &a\end{smallmatrix}\right)
\qquad
T = \left(\begin{smallmatrix} a&c \\ b &d\end{smallmatrix}\right)
\qquad
T^{-1} = \left(\begin{smallmatrix} d&-c \\ -b &a\end{smallmatrix}\right)
\end{equation}
We need only prove that $S(Q_1) \cap T(Q_1) = \emptyset$.
Since $Q_3=-Q_1$, this immediately yields $S(Q_3) \cap T(Q_3) = \emptyset$.
Consider the matrix $A=\left(\begin{smallmatrix} 0&1 \\ -1 &0\end{smallmatrix}\right)$
that maps $Q_1$ bijectively to $Q_2$.
Then
$$|S^{-1}(Q_2) \cap T^{-1}(Q_2)| = |A^{-1} S^{-1} A(Q_1) \cap A^{-1} T^{-1} A (Q_1)| = |T(Q_1) \cap S(Q_1)| = 0\,.$$
Similarly, since $Q_2=-Q_4$, this yields $S^{-1}(Q_4) \cap T^{-1}(Q_4) = \emptyset$ as well.

\medskip

Now suppose that $S(Q_1) \cap T(Q_1) \neq \emptyset$.  We will derive a contradiction.
Restating our assumption, there exists $(x,y) \in Q_1$ with $S^{-1} T(x,y) \in Q_1$.
This implies that
\begin{eqnarray}
(ad-b^2)x + d(c-b)y & >& 0 \label{eq:con1}\\
a(b-c)x + (ad-c^2)y &\geq & 0 \,. \label{eq:con2}
\end{eqnarray}

Note that $b \neq c$ since, by assumption, $S^{\top} \neq S$.
Also, $ad \neq 0$, since in this case $bc=-1$, which is impossible under our assumption that $b,c \geq 0$.

If $ad=c^2$ then $1=ad-bc=c(c-b)$ which implies that $c=1$ and $b=0$.  This yields $-ax \geq 0$ in \eqref{eq:con2},
which is impossible since $(x,y) \in Q_1 \implies x > 0$.

If $ad=b^2$ then $1=ad-bc=b(b-c)$, which implies that $c=0$ and $b=1$.
Altogether, in this case, we have $S = \left(\begin{smallmatrix} 1&1 \\ 0 &1\end{smallmatrix}\right)$.
Here we can conclude that $S(Q_1) \cap T(Q_1) = \emptyset$ because
$S$ maps points of $Q_1$ strictly below the line $y=x$ and $T$ maps points
of $Q_1$ above (or onto) the line $y=x$.

To summarize, we are left to deal with the case
\[
b \neq c, \ \ a > 0,\ \  d > 0,\ \  ad \neq b^2,\ \  ad \neq c^2\,.
\]

If $b > c$ then $ad-b^2 < ad-bc = 1$ which implies $ad-b^2 < 0$ since $ad \neq b^2$.
In this case, $d(c-b) < 0$ as well.  Thus if \eqref{eq:con1} holds, then $x=y=0$.
Similarly, if $c > b$, then $ad-c^2 < ad-bc = 1$ hence $ad-c^2 < 0$ and $a(b-c) < 0$,
implying $x=y=0$.  We conclude that $S(Q_1) \cap T(Q_1) = \emptyset$.
\end{proof}

\begin{corollary}\label{cor:positive}
If $S \in GL_2(\Z)$ has all non-negative entries, $S \neq S^{\top}$, and $\det(S)=1$, then
for any subset $A \subseteq \Z^2 \setminus \{0\}$,
$$
|S(A) \cup S^{\top}(A) \cup S^{-1}(A) \cup S^{-\top}(A)| \geq 2|A|\,.
$$
In particular, $h(G^S) > 0$.
\end{corollary}

\begin{proof}
In this case, we have $S(Q_1), S^{\top}(Q_1) \subseteq Q_1$, $S(Q_3), S^{\top}(Q_3) \subseteq Q_3$,
$S^{-1}(Q_2), S^{-\top}(Q_2) \subseteq Q_2$, and $S^{-1}(Q_4), S^{-\top}(Q_4) \subseteq Q_4$.
Thus \pref{lem:transpose} yields the desired result.
\end{proof}

To handle the case of general $S \in GL_2(\Z)$, it will help to have the following well-known fact.

\begin{lemma}\label{lem:simulate}
Consider two infinite graphs $G=(V,E)$ and $G'=(V,E')$ on the same countable index set $V$,
both of which have uniformly bounded degree.
Suppose there is a number $k \in \vvmathbb{N}$ such that
that for every $\{x,y\} \in E$, there is a path of length at most $k$ between $x$ and $y$ in $G'$.
Then $h(G) > 0$ implies $h(G') > 0$.
\end{lemma}

\begin{proof}
Let $\Delta$ be a uniform upper bound on the degree of vertices in $G$ and $G'$.
For a subset $U \subseteq V$ and $j \geq 1$, write $N^j_{G'}(U) \subseteq V$ for the set of vertices
within distance $j$ of the set $U$ in $G'$.

Now, suppose that $h(G')=0$.  In that case, for every $\e > 0$, there exists
a finite subset $U \subseteq V$ such that $|N_{G'}^1(U)| \leq (1+\e) |U|$.
In particular, this implies that $|N_{G'}^k(U)| \leq (1+\e \Delta^k)|U|$.  But, by our assumptions on $G$ and $G'$, this
implies
$$
|E(U, \bar U)| \leq \Delta (|N^k_{G'}(U)| - |U|) \leq \e \Delta^{k+1} |U|\,.
$$
Letting $\eps \to 0$ shows that $h(G)=0$ as well.
\end{proof}

The following two simple lemmas give conditions under which $G^{S,T}$ has Cheeger constant zero.

\begin{lemma}\label{lem:one}
For any $S \in GL_2(\Z)$, we have $h(G^{S,S^{-1}})= h(G^{S,-S^{-1}})=0$.
\end{lemma}

\begin{proof}
Let $G=G^{S,\pm S^{-1}}$ have edge set $E$.
Consider the sets
$\{U_k \subseteq \Z^2\}$ given by $$U_k = \{(j,0), S(j,0), \ldots, S^k(j,0) : j \in \{-1,1\}\}\,.$$
If $\sup_k |U_k| < \infty$, then clearly $h_{G^S}(U_k)=0$ for some $k$.  Otherwise,
since $|E(U_k)| \leq 4$, it must be that
$h_{G^S}(U_k) \to 0$ as $k \to \infty$, implying that $h(G^S)=0$.
\end{proof}

\begin{lemma}\label{lem:fourth}
Suppose $S,T \in GL_2(\Z)$ satisfy $S^4=T^4=I$.  Then $h(G^{S,T})=0$.
\end{lemma}

\begin{proof}
First, an elementary calculation shows that if $A \in GL_2(\Z)$ satisfies $\det(A)=1$ and $A^2=I$, then $A \in \{-I,I\}$.
Thus $S^2,T^2 \in \{-I,I\}$.
So for any $j_1, k_1, j_2, k_2, \ldots, j_m, k_m \in \Z$, we have
$$S^{j_1} T^{k_1} S^{j_2} T^{k_2} \cdots S^{j_m} T^{k_m} = (-1)^{i_0} T^{j_0} (ST)^j S^{k_0}\,.$$
for some $i_0,j_0,k_0 \in \{0,1\}$ and $j \in \vvmathbb{N} \cup \{0\}$.  Consider now the sets
$$
U_k = \left\{ (-1)^{i_0} T^{j_0} (ST)^j S^{k_0} (1,0) : i_0,j_0,k_0 \in \{0,1\} \textrm{ and } 0 \leq j \leq k \right\}\,.
$$
Letting $E^{S,T}$ denote the edge set of $G^{S,T}$, we have $|E^{S,T}(U_k, \bar U_k)| \leq 2 \cdot 8$ for every $k \geq 1$, and thus $h(G^{S,T}) = 0$.
\end{proof}

Finally, we complete the proof of \pref{thm:GS}.

\begin{proof}[Proof of \pref{thm:GS}]
Suppose that $S = \left(\begin{smallmatrix} a & b \\ c & d \end{smallmatrix}\right)\in GL_2(\Z)$ satisfies $S \neq S^{\top}$ and $\mathrm{tr}(S) \neq 0$, i.e.
$b \neq c$ and $a+d \neq 0$.  Let $T = S^{\top}$.
If $S$ has all non-negative or all non-positive entries, then the matrix
$S^2 = \left(\begin{smallmatrix} a^2+bc & b(a+d) \\ c(a+d) & bc+d^2 \end{smallmatrix}\right)$
has all non-negative entries, $\det(S^2)=1$, and $S^2 \neq (S^2)^{\top}$ by our initial assumptions.
Therefore by \pref{cor:positive}, we have $h(G^{S^2}) > 0$.
Now \pref{lem:simulate} implies $h(G^S) > 0$ as well.

If $ad > 0$ then $|\det(S)|=1$ implies $bc \geq 0$.  In this case, $S^{-1}$ has all non-negative
or all non-positive entries, hence $h(G^S) = h(G^{S^{-1}}) > 0$ by the preceding paragraph.

Thus we are left to deal with the case $ad \leq 0$.
But now consider the matrix $ST^{-1} =  \det(S) \left(\begin{smallmatrix} ad-b^2 & a(b-c) \\ d(c-b) & ad-c^2 \end{smallmatrix}\right)$.
We have $\det(S T^{-1})=1$ and $ST^{-1} \neq (ST^{-1})^{\top}$, by our initial assumptions that $b \neq c$ and $a+d \neq 0$.  Furthermore, the diagonal entries
of $ST^{-1}$ have the same sign, so our previous considerations yield $h(G^{ST^{-1}}) > 0$.
By \pref{lem:simulate}, this yields $h(G^S) > 0$ as well.

\medskip

To finish the proof, we must now show that if $S$ satisfies $S = S^{\top}$ or $\mathrm{tr}(S)=0$ then $h(G^S)=0$.
In the former case, we can apply \pref{lem:one}.
If $\mathrm{tr}(S)=0$, then $S^2 = \left(\begin{smallmatrix} a^2+bc & 0 \\ 0 & bc+d^2 \end{smallmatrix}\right) = \pm I$.
Similarly, $T^2 = \pm I$.
Thus $h(G^S)=0$ by \pref{lem:fourth}.
\end{proof}

\subsection{Conjugating by a reflection}
\label{sec:cai}

To further exhibit the flexibility of our method, we analyze the expansion
a different family of operators considered earlier by Cai \cite{Cai03}.
Let $R=\left(\begin{smallmatrix} 0&1 \\1 &0\end{smallmatrix}\right)$ and
for every $S \in GL_2(\Z)$, consider the graph
$$G^{S,RSR} =  \left(\Z^2 \setminus \{0\}, \left\{ \vphantom{\bigoplus}\{z,S z\}, \{z, RSR z\} : z \in \Z^2 \setminus \{0\}\right\}\right)\,.$$
Our goal is to prove the following analog of \pref{thm:GS}.

\begin{theorem}\label{thm:HS}
For any $S = \left(\begin{smallmatrix} a & b \\ c & d \end{smallmatrix}\right) \in GL_2(\Z)$, we have $h(G^{S,RSR}) > 0$ if and only if
$(a+d)(b+c) \neq 0$.
\end{theorem}

The next result follows from the preceding theorem and \pref{thm:maingen}

\begin{theorem}
For any $S = \left(\begin{smallmatrix} a & b \\ c & d \end{smallmatrix}\right) \in GL_2(\Z)$, $\{G_n^{S,RSR}\}$ is an expander
family if and only if $(a+d)(b+c) \neq 0$.
\end{theorem}

\begin{proof}
By \pref{thm:HS}, we have $h(G^{S^{\top}, RS^{\top} R}) > 0$.  Now \pref{thm:maingen} implies
that $\{G_n^{S,RSR}\}$ is
an expander family, noting that $(RSR)^{\top} = R S^{\top} R$.

On the other hand, suppose that $a+d=0$.  Then $S^4 = I = RS^4R$ so \pref{lem:fourthfin}
implies that $\{G_n^{S,RSR}\}$ is not an expander family.  If $b+c=0$ then $ST \in \{-I,I\}$,
so \pref{lem:onefin} implies the same.
\end{proof}

To illustrate another method of expansion analysis, we recall the following result of \cite{Cai03}.
Gabber and Galil \cite{GG81} proved this for $S = \left(\begin{smallmatrix} 1& 1 \\ 0 & 1 \end{smallmatrix}\right)$.

\begin{theorem}\label{thm:cai}
Consider any $S = \left(\begin{smallmatrix} a & b \\ c & d \end{smallmatrix}\right) \in GL_2(\Z)$ such that
$\det(S)=1$ and $|a+d|, |b+c| \geq 2$ are satisfied.  Then for any $z \in \Z^2 \setminus \{0\}$, one
of the following two conclusions holds for the set $$\left\{\|S z\|_{\infty}, \|S^{-1} z\|_{\infty}, \|RSR z\|_{\infty}, \|RS^{-1}R z\|_{\infty}\right\}\,.$$
Either three of the elements are strictly greater than $\|z\|_{\infty}$ or
at most two are equal to $\|z\|_{\infty}$ and the rest are strictly greater than $\|z\|_{\infty}$.
\end{theorem}

This rather immediately yields a positive Cheeger constant for $G^{S,RSR}$.

\begin{theorem}\label{thm:RSR}
Suppose that $S$ satisfies the assumptions of \pref{thm:cai}.  Then $h(G^{S,RSR}) > 0$.
\end{theorem}

\begin{proof}
For an edge $\{x,y\} \in E^{S,RSR}$, let
$$
\Delta(x,y) = \begin{cases}
0 & \|x\|_{\infty} = \|y\|_{\infty} \\
1 & \|x\|_{\infty} > \|y\|_{\infty} \\
-1 & \textrm{otherwise.}
\end{cases}
$$
Consider a finite set $U \subseteq \Z^2 \setminus \{0\}$.
Then by \pref{thm:cai},
$$
\sum_{x \in U} \sum_{ A \in \{S,RSR,S^{-1},RS^{-1}R\} } \Delta(x,Ax) \geq 2 |U|\,.
$$
On the other hand, whenever $x$ and $Ax$ are both in $U$, the total contribution
from the terms $\Delta(x,Ax)$ and $\Delta(Ax,x)$ is zero.  Thus at least $|U|/2$ elements
of $U$ have a neighbor outside $U$.  This implies that $h(G^{S,RSR}) > 0$.
\end{proof}

\begin{remark}
We observe that \pref{thm:cai} appears to be a genuinely different reason for expansion,
as an analysis akin to \pref{lem:transpose} does not appear to work in this setting when $ad \leq 0$.
To illustrate this, suppose that $S=\left(\begin{smallmatrix} a & b \\ c & d \end{smallmatrix}\right)$
and $a,b > 0$ and $c,d < 0$.  Setting $T=RSR$, one has $S(Q_1) \subseteq Q_2, S(Q_3) \subseteq Q_4,
S^{-1}(Q_1) \subseteq Q_4, S^{-1}(Q_3) \subseteq Q_2, T(Q_1) \subseteq Q_4, T(Q_3) \subseteq Q_2,
T^{-1}(Q_1) \subseteq Q_2, T^{-1}(Q_3) \subseteq Q_4$.  Notice that unlike in the case of $T=S^{\top}$,
one can only restrict the images to a single quadrant when the domain is $Q_1$ or $Q_3$.  This seems
to elude the simple counting argument of \pref{lem:transpose} and \pref{cor:positive}.
\end{remark}

We can now prove our main theorem.

\begin{proof}[Proof of \pref{thm:HS}]
Suppose first that $S = \left(\begin{smallmatrix} a&b \\ c &d\end{smallmatrix}\right)\in GL_2(\Z)$ satisfies $\det(S)=1$ and $(a+d)(b+c) \neq 0$.
Consider the matrix $S(RSR) = \left(\begin{smallmatrix} b^2+ad & a(b+c) \\ d(b+c) & c^2+ad\end{smallmatrix}\right)$.
First, we have
\begin{equation}\label{eq:trace}
\mathrm{tr}(SRSR) =b^2+c^2+2ad = b^2 +c^2 +2(1+bc) = (b+c)^2 + 2 > 2\,,
\end{equation}
where we have used $ad-bc=1$.  Let $\left(\begin{smallmatrix} u&v \\ w &x\end{smallmatrix}\right)$ denote $SRSR$ and
note that \eqref{eq:trace} gives $u+x > 2$.

We have $(SRSR)^2 =  \left(\begin{smallmatrix} u^2 + vw & v(u+x) \\ w(u+x) & x^2 + vw\end{smallmatrix}\right).$
The sum of the diagonal entries of this matrix is
$$
u^2 + x^2 + 2vw = u^2 + x^2 + 2(ux - 1) \geq (u+x)^2 - 2 > 2\,,
$$
where we have used $1=\det(SRSR)=ux-vw$ and $u+x \geq 2$.  Furthermore,
the sum of the off-diagonal entries satisfies
$$
|(w+v)(u+x)| \geq 2|w+v| = 2|(a+d)(b+c)| \geq 2\,.
$$
where we have used the assumption that $(a+d)(b+c) \neq 0$.
Thus we can apply \pref{thm:RSR} to $(SRSR)^2$ to conclude that
$h(G^{(SRSR)^2, R(SRSR)^2R}) > 0$.
Noting that $$R(SRSR)^2R = R(SRSR)(SRSR)R = (RSR)S(RSR)S\,$$
we can apply \pref{lem:simulate} to conclude that $h(G^{S,RSR}) > 0$ as well.

Finally, consider the case $\det(S)=-1$ and $(a+d)(b+c) \neq 0$.  The
matrix $S^2 =  \left(\begin{smallmatrix} a^2+bc& b(a+d) \\ c(a+d) &d^2+bc\end{smallmatrix}\right)$
satisfies $\det(S^2)=1$.  The sum of the off-diagonal entries is $(b+c)(a+d) \neq 0$.  The
sum of the diagonal entries is $a^2+d^2+2bc = a^2+d^2+2(ad-1) = (a+d)^2 -2 \neq 0$.  Thus
the preceding paragraph implies that $h(G^{S^2,RS^2R}) > 0$.
Now \pref{lem:simulate} yields $h(G^{S,RSR}) > 0$ as well.

\medskip

We now address the cases where the Cheeger constant is zero.
Write $T=RSR$.
If $a+d=0$ then $S^2 = \pm I$ and $T^2 = \pm I$, so \pref{lem:fourth} yields $h(G^{S,T})=0$.
If $b+c=0$ then $ST = \left(\begin{smallmatrix} b^2+ad&0 \\ 0 & b^2+ad \end{smallmatrix}\right) = \pm I$,
so \pref{lem:one} yields $h(G^{S,T})=0$.
\end{proof}

\subsection{Transformations for which $\{G_n^{S,T}\}$ is not an expander family}

Here, we argue that if $T=S^{-1}$ or $S^4=T^4=I$, then the graphs $\{G_n^{S,T}\}$ do not form
expander families.  The arguments are related to \pref{lem:one} and \pref{lem:fourth},
respectively, but we must also address the isoperimetric properties of boxes under linear transformations.
To this end, we define for $L \geq 0$ the box $B_L = \{ (x,y) \in \vvmathbb{R}^2 : -L \leq x \leq L, -L \leq y \leq L \}$.
For a subset $\Omega \subseteq \vvmathbb{R}^2$, we write $[\Omega] = \Omega \cap \Z^2$.
We also use $E_{\Z^2}$ to denote the edge set of the canonical graph on the integer lattice
where $x,y \in \Z^2$ are connected by an edge if and only if $\|x-y\|_1=1$.
The next lemma follows from elementary geometric considerations.

\begin{lemma}\label{lem:euclid}
For every $S \in GL_2(\Z)$, there is a constant $c > 0$ such that the following holds.  For every $L \geq 0$,
$S(B_L)$ is a parallelogram with area $4L^2$ and perimeter at most $c L$.
Furthermore,
we have $\liminf_{L \to \infty} [S(B_L)]/L^2 > 0$ and $\limsup_{L \to \infty} |E_{\Z^2}([S(B_L)])|/L \leq c$.
\end{lemma}

We also have the following basic classification of matrices in $GL_2(\Z)$; see, e.g. \cite[Ch. 1]{Gunning62}.

\begin{lemma}\label{lem:class}
Every $S \in GL_2(\Z)$ satisfies exactly one of the following.
\begin{enumerate}
\item $S$ has order dividing 12.
\item $S$ is conjugate in $GL_2(\vvmathbb{R})$ to $\left(\begin{smallmatrix} \alpha & 0 \\  0 & \alpha^{-1} \end{smallmatrix}\right)$
for some $\alpha \in \vvmathbb{R}$ with $|\alpha|,|\alpha^{-1}| \neq 1$.
\item $S$ is conjugate in $GL_2(\vvmathbb{R})$ to $\pm 1\left(\begin{smallmatrix}  1 & \gamma \\ 0 &  1\end{smallmatrix}\right)$ for some $\gamma \in \vvmathbb{R}$.
\end{enumerate}
\end{lemma}

The next lemma demonstrates our approach to proving non-expansion.

\begin{lemma}\label{lem:onefin}
For any $S \in GL_2(\Z)$, if $T \in \{S^{-1}, -S^{-1}\}$, it holds that $\{G_n^{S,T}\}$ is not an expander family.
\end{lemma}

\begin{proof}
For $T \in \{S^{-1},-S^{-1}\}$,
let $\bar G$ have vertex set $\Z^2$ and edge set $E = E^{S,T} \cup E_{\Z^2}$.
We will prove that $h(\bar G)=0$.  This is sufficient to show that $\{G_n^{S,T}\}$ is not an expander family.
Indeed, if $\{U_k\}$ is a sequence of finite sets with $h_{\bar G}(U_k) \to 0$, then for each $k$
one can choose the modulus $n$ large enough to avoid ``wrap around,'' yielding $h_{G_n^{S,T}}(U_k)=h_{\bar G}(U_k)$,
where we consider $U_k$ as a set of vertices in $G_n^{S,T}$ by reducing modulo $n$.

For $k \in \vvmathbb{N}$ and $L \geq 0$, consider the sets
$\{U_k(L) \subseteq \Z^2\}$ given by $$U_k(L) = [B_L] \cup [S(B_L)] \cup [S^2(B_L)] \cup \cdots \cup [S^k(B_L)] \,.$$
Observe that $B_L=-B_L$.

If we are in case (i) of \pref{lem:class}, then $U_{k_0}=U_{k_0}+1$ for some finite $k_0$.
So by \pref{lem:euclid}, we have $\liminf_{L \to \infty} |U_{k_0}(L)| \geq 4 L^2$, while $E^{S,T}(U_{k_0}(L)) = \emptyset$
and $\limsup_{L \to \infty} |E_{\Z^2}(U_{k_0}(L))| \leq c L$, where $c$ is some constant
depending on $S$ and $k_0$.   Thus $\lim_{L \to \infty} |E(U_{k_0}(L))|/|U_{k_0}(L)| = 0$ and $h(\bar G)=0$.

Now suppose that we are in case (ii) of \pref{lem:class} and, without loss of generality, $|\alpha| > 1$.  In this case, for some constant $\e > 0$ (depending possibly on $S$)
and every $k \in \vvmathbb{N}$,
we have \begin{equation}\label{eq:voleq}
\liminf_{L \to \infty} |U_k(L)|/L^2 \geq \e k\,.
\end{equation}
This follows
because the eccentricity of the parallelogram $S^k (B_L)$ grows exponentially fast; in fact, proportional to $|\alpha|^k$.
Similarly, in case (iii) of \pref{lem:class}, there is an $\e > 0$ (depending on both $S$)
such that $\liminf_{L \to \infty} |U_k(L)|/L^2 \geq \e k$.
To see this, it suffices to consider the case $\gamma=1$ in (iii) (since $\e$ can be depend on $\gamma)$.
In that case, the set $A_k = B_L \cup S(B_L) \cup \cdots \cup S^k(B_L)$ contains
an isosceles triangle whose corners are $\{(0,0), (kL,L), (-kL,L)\}$,
thus the volume of $A_k$ is at least $kL^2$.
Therefore \eqref{eq:voleq} again holds.

On the other hand, from \pref{lem:euclid} it follows that for some constant $c > 0$ (depending on $S$ and $k$),
$\limsup_{L \to \infty} |E_{\Z^2}(U_k(L))|/L \leq c$
and $\limsup_{L \to \infty} |E^{S,T}(U_k(L))|/L^2 \leq c$.
Therefore,
$$
\limsup_{L \to \infty} \frac{|E(U_k(L))|}{|U_k(L)|} \leq \frac{c}{\e k}\,.
$$
Taking $k \to \infty$ shows that $h(\bar G)=0$.

Finally, suppose that $S$ satisfies case (iii) of \pref{lem:class}.
\end{proof}

\begin{lemma}\label{lem:fourthfin}
Suppose $S,T \in GL_2(\Z)$ satisfy $S^4=T^4=I$.  Then $\{G_n^{S,T}\}$ is not an expander family.
\end{lemma}

\begin{proof}
Let $\bar G$ have vertex set $\Z^2$ and edge set $E = E^{S,T} \cup E_{\Z^2}$.
As in \pref{lem:onefin}, it will suffice to show that $h(\bar G)=0$.

As in \pref{lem:fourth},
an elementary calculation shows that if $A \in GL_2(\Z)$ satisfies $\det(A)=1$ and $A^2=I$, then $A \in \{-I,I\}$.
Thus $S^2,T^2 \in \{-I,I\}$.
So for any $j_1, k_1, j_2, k_2, \ldots, j_m, k_m \in \Z$, we have
$$S^{j_1} T^{k_1} S^{j_2} T^{k_2} \cdots S^{j_m} T^{k_m} = (-1)^{i_0} T^{j_0} (ST)^j S^{k_0}\,.$$
for some $i_0,j_0,k_0 \in \{0,1\}$ and $j \in \vvmathbb{N}$.  Consider now the sets
$$
U_k(L) = \left\{ [T^{j_0} (ST)^j S^{k_0} B_L] : j_0,k_0 \in \{0,1\} \textrm{ and } 0 \leq j \leq k \right\}\,.
$$

We can apply \pref{lem:class} to the matrix $ST$;
the resulting case analysis is essentially the same
as \pref{lem:onefin}.
\end{proof}

\subsection*{Acknowledgements}

I am grateful to Mohammad Moharrami, Yuval Peres, and
the audience at Microsoft Research for many helpful comments.

\bibliographystyle{alpha}
\bibliography{margulis}

\end{document}